\documentclass[15pt]{amsart}
 
\usepackage[english]{babel}
\usepackage[all]{xypic}
\usepackage{amsmath}
\usepackage{amsthm}
\usepackage{amsfonts}
\usepackage{hyperref}
\usepackage{amssymb}
\usepackage{stmaryrd}
\usepackage{enumerate}

\newcommand{\CC}{\mathbb{C}}

\newcommand{\NN}{\mathbb{N}}

\newcommand{\RC}{\mathcal{R}}

\newcommand{\CCC}{\mathcal{C}}
\newcommand{\DC}{\mathcal{D}}
\newcommand{\NC}{\mathcal{N}}
\newcommand{\AC}{\mathcal{A}}
\newcommand{\BC}{\mathcal{B}}

\providecommand{\tor}{\operatorname{Tor}}
\providecommand{\betti}{b^{(2)}}
\providecommand{\boldP}{\mathbf{P}}
\providecommand{\boldT}{\mathbf{T}}

\providecommand{\gr}{{\underline{G}}}

\providecommand{\lgrb}{{L^\infty (\gr^0)}}
\providecommand{\lxb}{{L^\infty (X)}}

\providecommand{\lxbg}{{\lxb\mspace{-2mu}\ast\mspace{-2mu} G}}

\providecommand{\curver}{\mspace{-4mu}\curvearrowright\mspace{-4mu}}

\providecommand{\trace}{\operatorname{tr}}

\providecommand{\im}{\operatorname{im}}

\providecommand{\colim}{\operatorname{colim}}
\providecommand{\coker}{\operatorname{coker}}

\theoremstyle{definition}
\newtheorem{defi}{Definition}[section]

\newtheorem{rem}[defi]{Remark}

\theoremstyle{plain}
\newtheorem{lem}[defi]{Lemma}
\newtheorem{thm}[defi]{Theorem}

\begin{document}
\title{L$^2$-Betti Numbers of Discrete Measured Groupoids} 
\author{Roman Sauer}
\email{roman.sauer@uni-muenster.de}
\urladdr{wwwmath.uni-muenster.de/math/inst/reine/inst/lueck/homepages/roman\_sauer/sauer.html}     
\address{FB Mathematik, Universit\"at M\"unster, Einsteinstr. 62,
  48149 M\"unster\\Germany} 

\begin{abstract} 
There are notions of L$^2$-Betti numbers for discrete groups
(Cheeger-Gromov, L\"uck), for 
type $II_1$-factors (recent work of Connes-Shlyakhtenko) and 
for countable standard equivalence relations (Gaboriau). Whereas the 
first two are algebraically defined using L\"uck's dimension theory, 
Gaboriau's definition of 
the latter is inspired by the work of Cheeger and Gromov. In this work 
we give a definition of L$^2$-Betti numbers of discrete measured
groupoids that is based on L\"uck's
dimension theory, thereby encompassing the cases of groups,
equivalence relations and holonomy groupoids with an invariant measure
for a complete transversal. 
We show that with our definition, like with Gaboriau's, 
the L$^2$-Betti numbers $\betti_n(G)$ of a
countable group $G$ coincide with the L$^2$-Betti numbers
$\betti_n(\RC)$ of the orbit equivalence relation $\RC$ of 
a free action of $G$ on a probability space. This yields a new proof
of the fact the $L^2$-Betti numbers of groups with 
orbit equivalent actions coincide. 
\end{abstract}

\maketitle
\section{Introduction and Statement of
  Results}\label{sec:intr-stat-results}
In~\cite{cheeger-gromov(1986)} Cheeger and Gromov defined $L^2$-Betti
numbers for arbitrary countable discrete groups. In a series of
papers~\cite{Lueck(1997a)},~\cite{Lueck(1998a)},~\cite{Lueck(1998b)} 
L\"uck put the theory of $L^2$-Betti numbers into a completely algebraic
framework by introducing a dimension function $\dim_{\NC(G)}$ for
arbitrary modules over the group von Neumann algebra $\NC(G)$: The 
$L^2$-Betti numbers $\betti_n(G)$ of a group $G$ can be read off from
the group homology as 
\[\betti_n(G)=\dim_{\NC(G)}\bigl (\tor_n^{\CC G}(\NC(G), \CC)\bigr
).\]
More recently,
in an influential paper of Gaboriau~\cite{gaboriau(2002b)} the notion 
of $L^2$-Betti numbers for countable standard equivalence relations  
was introduced. Their construction is motivated by the one of 
Cheeger and Gromov. This article provides an homological-algebraic 
definition of $L^2$-Betti numbers for countable standard equivalence
relations and, more generally, for discrete measured groupoids. 
In section~\ref{sec:discr-meas-group} 
we will introduce algebraic objects for a discrete measured groupoid
$\gr$ like the groupoid ring $\CC\gr$ 
which are analogous to the group case, and we define 
the $L^2$-Betti numbers $\betti_n(\gr)$ of $\gr$ in~\ref{betti-definition} as 
\[\betti_n(\gr)=\dim_{\NC(\gr)}\bigl (\tor_n^{\CC\gr}(\NC(\gr),
L^\infty(\gr^0))\bigr ).\]

Here $\gr^0\subset\gr$ denotes the subset of objects in the groupoid. 
Denote the source and target maps by $s$ resp.~$t$. We show
in~\ref{betti-restriction} the following formula for the restriction. 
\begin{thm}\label{einl1}
Let $\gr$ be a discrete measured groupoid, and let $A\subset\gr^0$ 
be a Borel subset such that $t(s^{-1}(A))$ has full measure in 
$\gr^0$. Then 
\[\betti_n(\gr)=\mu (A)\cdot\betti_n (\gr_{\vert A}).\]
\end{thm}

An essentially free 
action of a countable group $G$ on a standard Borel probability space
$X$ by 
measure preserving Borel automorphisms is called a 
\emph{standard action}, and $X$ is called a \emph{probability $G$-space}. 
Let $\RC$ be the orbit equivalence relation 
of a standard action of the group $G$. 
Then $\RC$ is an example of a discrete 
measured groupoid. The following theorem is proved
in~\ref{betti-mainthm} by homological-algebraic 
methods, which might be
useful in other contexts dealing with homological algebra of finite
von Neumann algebras. 
\begin{thm}\label{einl2} $\betti_n(G)=\betti_n(\RC)$.\end{thm}
The standard actions $G\curver X$, $H\curver Y$ are called 
\emph{orbit equivalent} if there exists a measure-preserving 
Borel isomorphism $f:X\rightarrow Y$ that maps orbits bijectively onto 
orbits. They are called \emph{weakly orbit equivalent} if there are 
Borel subsets $A\subset X$, $B\subset Y$ (equipped with the normalized
measures) meeting almost every orbit
and a measure-preserving Borel isomorphism $f:A\rightarrow B$ that
maps orbits bijectively onto orbits. 
The theorems above immediately imply  
\begin{thm}[\cite{gaboriau(2002b)}]\label{einl3} Assume 
$G$ and $H$ possess weakly orbit equivalent standard actions. 
Then there is a $C>0$ such that $\betti_n(G)=C\cdot\betti_n(H)$ holds for 
$n\ge 0$. Further, if the actions are orbit equivalent then 
$\betti_n(G)=\betti_n(H)$ for $n\ge 0$. 
\end{thm}

As a remark, two countable groups possess weakly orbit equivalent 
standard actions if and only if they are measure
equivalent~\cite{Furman(1999a)},\cite{Furman(1999b)}. 
In~\cite{gaboriau(2002b)} Gaboriau proves
theorems~\ref{einl1},~\ref{einl2},~\ref{einl3} for his version 
of $L^2$-Betti numbers of countable equivalence relations. 
We remark that our methods are independent of
Gaboriau's results. In particular, we do \textit{not} show 
the equality of our definition and Gaboriau's one 
to deduce~\ref{einl2}. Of course, as a consequence of~\ref{einl2} Gaboriau's
definition and ours 
coincide for orbit equivalence
relations of free measure-preserving group actions 
and presumably coincide for every countable standard equivalence
relation. 

Popa showed 
that a type $II_1$ factor $\BC$ satisfying some rigidity and 
compact approximation properties has only one Cartan subalgebra 
$\AC$~\cite{popa(2002)}. Thus the $L^2$-Betti numbers of the countable standard
equivalence relation associated to the inclusion $\AC\subset\BC$ 
are isomorphism invariants of the factor $\BC$. 
Recently, Connes and Shlyakhtenko~\cite{connesshlyakhtenko(2003)} defined 
$L^2$-Betti numbers of an arbitrary type $II_1$-factor $\NC$ 
in a very different fashion as 
\[\betti_n(\NC)=\dim_{\NC\bar{\otimes}\NC^o}\bigl
(\tor_n^{\NC\otimes\NC^o} (\NC\bar{\otimes}\NC^o, \NC)\bigr ).\] 
The motivation for this article is twofold. 
First the definition of $L^2$-Betti numbers of  
discrete measured groupoids encompasses the examples 
of groups, countable standard equivalence relations and holonomy
groupoids with an invariant measure in a unifying
way. Further, using our definition of $\betti_n(\gr)$, 
one can put Popa's $L^2$-Betti numbers of $II_1$-factors
(if defined) into a homological algebra framework, which might be helpful 
in order to understand better the relation between these and the 
Connes-Shlyakhtenko $L^2$-Betti numbers. 

The results in this paper are part of the author's thesis. I want to 
thank my adviser Wolfgang L\"uck for constant support and
encouragement. 

\section{Review of Dimension Theory of finite von Neumann
  Algebras}\label{sec:revi-dimens-theory}
In this section we review the dimension function for arbitrary modules 
over a finite von Neumann algebra (modules in the algebraic
sense) and its basic properties. Further, we prove some additional 
properties concerning restrictions of von Neumann algebras and give a 
new criterion for a module to be zero dimensional.  
In the sequel let $\AC$ be a finite von Neumann algebra with 
normalized trace $\trace_\AC$. 
The dimension (function) $\dim_\AC (M)\in [0,\infty]$ of a module 
$M$ over $\AC$ was introduced by
W.~L\"uck in~\cite{Lueck(1998a)},~\cite{Lueck(1998b)}. 
For a finitely generated projective $\AC$-module $P$ 
choose an idempotent matrix $A=(A_{ij})\in M_n(\AC)$ such that 
$P\cong\AC^n\cdot A$. Then the dimension $\dim_\AC (P)$ is 
defined as 
\[\dim_\AC (P)=\sum_{i=1}^n\trace_\AC (A_{ii})\in [0,\infty).\]
For an arbitrary $\AC$-module $N$ the dimension is then defined as 
\[\dim_\AC (N)=\sup\{\dim_\AC (P);~\text{$P\subset N$ finitely
  generated projective submodule}\}\in [0,\infty].\]
For the following fundamental 
theorem see~\cite[Theorem 0.6]{Lueck(1998a)}, 
also~\cite[Theorem 6.7 on p.~239]{lueck(2002)}. 

\begin{thm}\label{betti-dim-properties}
The dimension function $\dim_\AC$ satisfies the following 
properties. 
\begin{enumerate}[(i)]
\item A projective $\AC$-module $P$ is trivial if and only 
if $\dim_\AC (P)=0$. 
\item Additivity.\\
If $0\rightarrow A\rightarrow B\rightarrow C\rightarrow 0$ 
is an exact sequence of $\AC$-modules then 
\[\dim_\AC (B)=\dim_\AC (A)+\dim_\AC (C)\]
holds, where we put $\infty+r=r+\infty=\infty$ for $r\in [0,\infty]$. 
\item Cofinality.\\
Let $M=\bigcup_{i\in I} M_i$ 
be a directed union of submodules $M_i\subset M$. Then 
\[\dim_\AC (M)=\sup_{i\in I}\{\dim_\AC (M_i)\}.\]
\end{enumerate}
\end{thm}

\begin{defi}An $\AC$-homomorphism $f:M\rightarrow N$ between $\AC$-modules 
$M$, $N$ is called a \emph{$\dim_\AC$-isomorphism} if 
$\dim_\AC (\ker f)=\dim_\AC (\coker f)=0$. 
\end{defi}

There is a suitable localization of the category 
of $\AC$-modules in which $\dim_\AC$-isomorphisms become 
isomorphisms. Let us recall the relevant notions.  

Let $\CCC$ be an abelian category. A non-empty full subcategory 
$\DC$ of $\CCC$ is called a \emph{Serre subcategory} if for all short exact 
sequences in $\CCC$ 
\[0\rightarrow M'\rightarrow M\rightarrow M''\rightarrow 0\] 
$M$ belongs to $\DC$ if and only if both $M'$ and $M''$ do. 
Then there is a quotient category $\CCC/\DC$ with the same objects as 
$\CCC$ and a functor $\pi:\CCC\rightarrow\CCC/\DC$. Moreover, 
$\CCC/\DC$ is abelian, $\pi$ is exact and $\pi(f)$ is an isomorphism 
if and only if $\ker (f)$ and $\coker (f)$ lie in $\DC$ for a morphism 
$f$ in $\CCC$. 
The properties in~\ref{betti-dim-properties} imply that the subcategory 
of zero-dimensional $\AC$-modules is a Serre subcategory. 
This has useful consequences. For instance, there is a 5-lemma for 
$\dim_\AC$-isomorphisms because there is one for general abelian 
categories. 

In~\cite[Lemma 3.4]{Lueck(1997a)} it is proved that 
any finitely generated $\AC$-module $N$ 
splits as $N=\boldP N\oplus\boldT N$ 
where $\boldP N$ is finitely generated projective and 
$\boldT N$ is the kernel of the canonical homomorphism 
$N\rightarrow N^{**}$ into the double dual, mapping 
$x\in N$ to $N^*\rightarrow\AC, f\mapsto f(x)$. 
Furthermore, $\dim_\AC (\boldT N)=0$ holds. 
The modules $\boldP N$ and $\boldT N$ are called the 
\emph{projective} resp.~\emph{torsion part of $N$}. 
Further, if $N$ is a finitely presented $\AC$-module then the 
torsion part 
$\boldT N$ possesses an exact resolution of 
the form 
\[0\rightarrow\AC^n\overset{r_A}{\longrightarrow}\AC^n\rightarrow\boldT N
\rightarrow 0,\]
where $r_A$ is right multiplication with a positive $A\in M_n(\AC)$. 

The next lemma is exercise 6.3 (with solution on p.~530) 
in~\cite[p.~289]{lueck(2002)}. 
It is formulated for group von Neumann algebras there, but the proof 
is exactly the same for arbitrary finite von Neumann algebras. 

\begin{lem}\label{betti-dim-exercise}
Let $M$ be a submodule of a finitely generated projective $\AC$-module $P$. 
For every $\epsilon>0$, 
there exists a submodule $P'\subset M$ that is a direct summand in $P$ 
and satisfies $\dim_\AC (M)\le\dim_\AC (P')+\epsilon$. 
\end{lem}

The following theorem is a local criterion for a module to be 
zero dimensional. 

\begin{thm}\label{betti-dim-localcriterion}
Let $M$ be an $\AC$-module. Its dimension  
$\dim_\AC (M)$ vanishes if and only if for every element $m\in M$ there is 
a sequence $p_i\in\AC$ of projections such that 
\[\lim_{i\rightarrow\infty}\trace_\AC (p_i)=1~\text{ and }~p_i\cdot m=0\text{
  for all $i\in\NN$.} \] 
Furthermore, if $q\in\AC$ is a given projection with $qm=0$ for an element 
$m$ in $M$ with $\dim_\AC (M)=0$, then 
the sequence $p_i$ can be chosen such that $q\le p_i$. 
\end{thm}

\begin{proof}
First assume $\dim_\AC (M)=0$. Consider an element $m\in M$, and let $q\in\AC$ 
be a projection such that $qm=0$. For a given $\epsilon>0$, 
we want to construct a projection $p\in\AC$ such that 
$\trace_\AC (p)\ge 1-\epsilon$, $p\cdot m=0$ and $p\ge q$. 
Let $\langle m\rangle\subset M$ be the submodule generated by $m$. 
We have the epimorphism
\[\phi:\AC(1-q)\rightarrow\langle m\rangle,~ n(1-q)\mapsto nm\]
and $\dim_\AC (\ker\phi)=\dim_\AC(\AC(1-q))-\dim_\AC(\langle
m\rangle)=1-\trace_\AC (q)$. 
By~\ref{betti-dim-exercise} there is a submodule $P\subset\ker\phi$ such 
that $P$ is a direct summand in $\AC(1-q)$ and $\dim_\AC(\ker\phi)\le\dim_\AC
(P) +\epsilon$. 
Hence $\AC q\oplus P\subset\AC q\oplus\AC (1-q)=\AC$ is a direct 
summand in $\AC$, i.e.~it is has the form $\AC p$ for a projection 
$p$. Its trace is 
$\trace_\AC (p)=\trace_\AC (q)+\dim_\AC (P)\ge 1-\epsilon$. 
Moreover, $\AC q\subset\AC p$ implies $qp=q$, i.e.~$q\le p$, and 
$pm=0$ is obvious. 
\smallskip\\
Now we prove the converse. 
It suffices to prove that $M$ has no non-trivial finitely generated 
projective submodules. Suppose $P\subset M$ is a non-trivial 
finitely generated projective submodule. Then there is a non-trivial 
$\AC$-homomorphism $\phi:P\rightarrow\AC$. Choose non-zero element 
$y=\phi(x)\ne 0$ in the image of $\phi$. 
There is a sequence of projections $p_i\in\AC$ such 
$\trace_\AC (p_i)\rightarrow 1$ and $p_i\cdot x=0$. In particular, 
$p_i\cdot y=\phi(p_i\cdot x)=0$ yielding $y=0$. Hence 
no such non-trivial $P$ can exist. 
\end{proof}

\begin{thm}\label{dim-functor}
Let $\AC$ and $\BC$ be finite von Neumann algebras, and let 
$F$ be an exact functor from the category of $\AC$-modules 
to the category of $\BC$-modules which preserves colimits. 
Assume there is a constant $C>0$ such that 
\begin{equation}\label{eq-forprojective}
\dim_\BC (F(P))=C\cdot \dim_\AC(P)
\end{equation}
holds for every finitely generated projective $\AC$-module $P$. Then 
$\dim_\BC (F(M))=C\cdot\dim_\AC (M)$ holds for every $\AC$-module 
$M$. 
\end{thm}

\begin{proof} 
We prove this for finitely presented, finitely generated and 
arbitrary modules, subsequently.\\
Step 1: Let $M$ be a finitely presented $\AC$-module. 
Then $M$ splits as $M=\boldP M\oplus\boldT M$, where $\boldP M$ is 
projective and $\boldT M$ admits an exact resolution 
\[0\rightarrow\AC^n\rightarrow\AC^n\rightarrow\boldT M\rightarrow 0.\]
Additivity yields $\dim_\AC (\boldT M)=0$. Applying the exact functor 
$F$ to this resolution, additivity also implies $\dim_\BC (F(\boldT M))=0$. 
We have $C\cdot\dim_\AC (\boldP M)=\dim_\BC (F(\boldP M))$ by
assumption. Hence, 
$C\cdot\dim_\AC (M)=\dim_\BC (F(M))$. \\
Step 2: Let $M$ be finitely generated. 
Then there is a finitely generated 
free $\AC$-module $P$ with an epimorphism $P\rightarrow M$. Let 
$K$ be its kernel. $K$ can be written as the directed union of its 
finitely generated submodules $K=\bigcup_{i\in I} K_i$. 
By cofinality and additivity (\ref{betti-dim-properties}) we conclude  
\begin{equation*}
\begin{split}
\dim_\AC (M) =\dim_\AC (P)-\dim_\AC (K)
           &=\dim_\AC (P)-\sup_{i\in I}\left \{\dim_\AC (K_i)\right\}\\
           &=\inf_{i\in I} \left\{\dim_\AC (P)-\dim_\AC (K_i)\right\}\\
           &=\inf_{i\in I} \left\{\dim_\AC (P/K_i)\right\}
\end{split}
\end{equation*}
We have $F(K)=\colim_{i\in I} F(K_i)$ with injective structure maps because 
$F$ is colimit-preserving and exact. Thus we can conclude similarly to 
obtain 
\[\dim_\BC (F(M))=\inf_{i\in I} \left\{\dim_\BC (F(P/K_i))\right\}.\]
Then the claim follows from the first step. \\
Step 3: Let $M$ be an arbitrary module.
Every module is the directed union of its finitely generated submodules, 
which reduces the claim to the preceding step due to cofinality. 
\end{proof}

The following theorem is proved in~\cite{Lueck(1998a)} for inclusions 
of group von Neumann algebras induced by group inclusions. Using the 
previous theorem it suffices to prove it for finitely generated
projective modules, which is easy. 

\begin{thm}\label{betti-dim-dimpreserving}
Let $\phi:\AC\rightarrow\BC$ be a trace-preserving $*$-homomorphism between 
finite von Neumann algebras $\AC$ and $\BC$. 
Then for every $\AC$-module $N$ we have  
\[\dim_\AC (N)=\dim_\BC \left (\BC\otimes_\AC N\right ).\]
\end{thm}

We recall some definitions and easy facts about 
Morita equivalence of rings. 
For details and proofs we refer to~\cite[section 18]{lam(1999)}. 
Two rings are called \emph{Morita equivalent} if there exists a 
category equivalence, called \emph{Morita equivalence}, 
between their module categories. Every Morita equivalence is exact and 
preserves projective modules. 
An idempotent $p$ in a ring $R$ is called \emph{full} if 
the additive subgroup in $R$ generated by the elements 
$rpr'$ with $r,r'\in R$, denoted by $RpR$, coincides with $R$. 
In this case $R$ and $pRp$ are 
Morita equivalent, and the mutual inverse category equivalences are 
given by tensoring with the bimodules $Rp$ resp.~$pR$. 

\begin{defi}\label{dim-full}
Let $\AC$ be a finite von Neumann algebra, and let $R$ be a ring 
containing $\AC$ as a subring. 
An idempotent $p\in R$ is called \emph{$\dim_\AC$-full} 
if the inclusion of $RpR$ into $R$ is a $\dim_\AC$-isomorphism of 
$\AC$-modules. 
\end{defi}

Note for the following that if $p$ is a projection in the finite 
von Neumann algebra $\AC$ then $p\AC p=\{pap;a\in\AC\}$ is again 
a finite von Neumann algebra equipped with the normalized trace 
$\trace_{p\AC p}(x)=\frac{1}{\trace_\AC(p)}\trace_\AC(x)$. 

\begin{thm}\label{betti-reducedxactness}
Let $p$ be a $\dim_\AC$-full projection in $\AC$. Then $\AC p$ is a 
right flat $p\AC p$-module. 
\end{thm}

\begin{proof}
Consider the image $\bar{1}$ of $1$ in the cokernel of the inclusion 
$\AC p\AC\subset\AC$. By assumption, the cokernel has dimension 
zero. By the local criterion~\ref{betti-dim-localcriterion} there 
is a sequence $(p_i)_{i\in\NN}$ of projections in $\AC$ such that 
$p_i\bar{1}=0$, i.e.~$p_i\in\AC p\AC$, $\trace_\AC (p_i)\rightarrow 1$ 
and $p_i\ge p$. From $p\le p_i$ we get $p=p_ipp_i\in p_i\AC p_i$. 
Furthermore, $p=p_ipp_i$ and $p_i\in\AC p\AC$ imply 
\[p_i\in \left (p_i\AC p_i\right )p\left (p_i\AC p_i\right ),\]
hence $p$ is a full idempotent in $p_i\AC p_i$, and the rings 
$p_i\AC p_i$ and $p\AC p$ are Morita equivalent. 
Thus the right $p\AC p$-module $p_i\AC p$ is projective. 
The $p\AC p$-homomorphism
\[\AC p\longrightarrow \prod_{i\in\NN} p_i\AC p,~n\mapsto (p_in)_{i\in\NN}\]
is injective. 
Now we use the fact that a von Neumann algebra is a semihereditary
ring (see~\ref{groupoids-rings-semihereditary}). 
Over a semihereditary ring the property of being flat 
is inherited to products and submodules
by~\ref{groupoids-rings-flatnesssemihereditary}. Therefore 
the product on the right is flat, and its 
submodule $\AC p$ is flat as a $p\AC p$-module. 
\end{proof}

\begin{thm}\label{betti-dim-projectiondim0}
Let $p$ be a $\dim_\AC$-full projection in $\AC$. 
For every $p\AC p$-module $M$ we have 
\[\dim_\AC\left (\AC p\otimes_{p\AC p} M\right )=\trace_\AC (p)\dim_{p\AC p}
(M).\]
\end{thm}

\begin{proof}
Due to the preceding theorem, 
the functor $\AC p\otimes_{p\AC p}-$ is exact. 
By~\ref{dim-functor} it suffices to check the claim for a finitely generated 
projective $p\AC p$-module $M$. 
Let $A\in M_n(p\AC p)$ be an idempotent matrix such that 
\[M\cong\left (p\AC p\right )^n\cdot A.\]
Then $\dim_{p\AC p} (M)=\sum_{i=1}^n\trace_{p\AC p} (A_{ii})$ by definition. 
The $\AC$-module $\AC p\otimes_{p\AC p} M$ is finitely generated 
projective, and we have 
\[\AC p\otimes_{p\AC p} M\cong\left (\AC p\right )^n\cdot
A=\AC^n\cdot A .\]
For the right equal sign note that $(1-p)A=0$. Hence we can conclude
\begin{equation*}
\begin{split}
\dim_\AC\left(\AC p\otimes_{p\AC p} M\right )=
\sum_{i=1}^n\trace_{\AC} (A_{ii})&=\trace_\AC (p)\cdot\sum_{i=1}^n\trace_{p\AC
  p} (A_{ii})\\&=\trace_\AC (p)\cdot\dim_{p\AC p} (M).
\end{split}
\end{equation*}
\end{proof}

\begin{thm}\label{almost-morita}
Let $R$ be a ring containing $\AC$ as a subring, and let 
$p$ be a $\dim_\AC$-full idempotent in $R$. Then 
\begin{equation}\label{eq-dimsurj}
\phi:R p\otimes_{pR p}pM\rightarrow M,~n\otimes m\mapsto nm 
\end{equation}
is a $\dim_\AC$-isomorphism for every $R$-module $M$. 
\end{thm}

\begin{proof}
First we show that $\phi$ is $\dim_\AC$-surjective. 
The local criterion~\ref{betti-dim-localcriterion} 
applied to the cokernel of the inclusion 
$ R p R\subset R$ provides a sequence $p_i$ of projections in $\AC$, which 
lie in $RpR$ and satisfy $\trace_\AC (p_i)\rightarrow 1$. 
But every element in the cokernel of $\phi$ is annihilated by 
the $p_i$, so this yields $\dim_\AC (\coker\phi)=0$, again due 
to the local criterion. 

Now we can prove that $\phi$ is a $\dim_\AC$-isomorphism. 
Consider the exact sequence 
\[0\rightarrow\ker\phi\rightarrow R p\otimes_{pRp}pM\rightarrow
M\rightarrow\coker\phi\rightarrow 0.\]
Applying the exact functor $pR\otimes_R -$ produces an isomorphism 
in the middle because $pR \otimes_R\left ( R p\otimes_{pR p}pM\right
)=pRp\otimes_{pRp} pM=pR\otimes_R M$.  
Hence $p\ker\phi=0$. 
Because $\phi$ is already shown to be $\dim_\AC$-surjective,  
we obtain $\dim_\AC (\ker\phi)=0$. Hence $\phi$ is a
$\dim_\AC$-isomorphism. 
\end{proof}

\begin{thm}\label{betti-dim-projectiondim}
Let $p$ be a $\dim_\AC$-full projection in $\AC$. 
Then for every $\AC$-module $M$ we have 
\[\dim_\AC (M)=\trace_\AC(p)\cdot\dim_{p\AC p} (pM).\]
\end{thm}

\begin{proof}
By the preceding theorem we have 
$\dim_\AC (\AC p\otimes_{p\AC p} pM)=\dim_\AC (M)$. Now the claim 
is obtained by~\ref{betti-dim-projectiondim0}. 
\end{proof}

\section{Discrete Measured Groupoids}\label{sec:discr-meas-group}

Discrete measured groupoids are generalizations of 
countable standard equivalence relations. The reference for 
the definitions and basic properties 
of the groupoid ring and von Neumann algebra associated to 
a countable standard equivalence relation is~\cite{feldman(1977b)}. 
We review the definitions in the more general setting of discrete 
measured groupoids.

In \cite[p.~82,88,89]{kechris(1995)} the following standard
measure-theoretical facts 
are proved: Any measurable subset of a standard Borel space is a standard 
Borel space. A bijective Borel map between standard Borel spaces is 
a Borel isomorphism. The image of an injective Borel map between
standard Borel spaces is Borel. We omit the proof of the 
following technical lemma. The case for countable standard equivalence 
relations is implicit in~\cite[prop.~2.3]{feldman(1977b)}. 
A detailed proof can be found in~\cite[theorem 1.3]{sauer(2002)}. 

\begin{lem}\label{groupoids-definition-kuratowski}
Let $f:X\rightarrow Y$ be a Borel map between standard Borel 
spaces such that the 
preimages $f^{-1}(\{y\})$, $y\in Y$, are countable. 
Then the image $f(X)$ is measurable in $Y$, and  
there is a countable partition of $X$ into measurable subsets 
$X_i$, $i\in\NN$, such that all $f_{|X_i}$ are injective, and 
$f_{|X_1}:X_1\rightarrow f(X)$ is a Borel isomorphism. 
If, in addition, there is some $N\in\NN$ such 
that $\#f^{-1}(\{y\})\le N$ for all $y\in Y$,  
then the partition can be chosen to have at most $N$ sets. 
\end{lem}

Recall that a \emph{groupoid} is a small category where all morphisms are 
invertible. We usually identify a groupoid $\gr$ with the set of 
its morphisms. The set of objects $\gr^0$ can be considered as 
a subset (via the identity morphisms). There are four canonical 
maps, namely 
\begin{enumerate}[]
\item the \emph{source
    map} $s:\gr\rightarrow\gr^0, (f:x\rightarrow y)\mapsto x$, 
\item the \emph{target map} $t:\gr\rightarrow\gr^0, (f:x\rightarrow y)\mapsto y$, 
\item the \emph{inverse map} $i:\gr\rightarrow\gr, f\rightarrow f^{-1}$ and 
\item the \emph{composition} 
$\circ: \gr^{(2)}:=\{(f,g)\in\gr\times\gr;~ s(f)=t(g)\}\rightarrow\gr,
  (f,g)\mapsto f\circ g$. 
\end{enumerate}
The composition will also be denoted by 
$g_1g_2$ instead of $g_2\circ g_1$. 
A \emph{discrete measurable groupoid} $\gr$ is a groupoid 
$\gr$ equipped with the structure of a standard Borel space 
such that the composition and the inverse map are Borel and 
$s^{-1}(\{x\})$ is countable for all $x\in\gr^0$. 
We remark that the source and target maps of a discrete measurable 
groupoid $\gr$ are measurable, $\gr^0\subset\gr$ is 
a Borel subset, and $t^{-1}(\{x\})$ is countable.  
Now let $\mu$ be a probability measure on 
the set of objects $\gr^0$ of a discrete measurable groupoid $\gr$. Then, 
for any measurable subset $A\subset\gr$, the function 
$\gr^0\rightarrow\CC,~x\mapsto \#\bigl (s^{-1}(x)\cap A\bigr )$ 
is measurable, and the measure $\mu_s$ on $\gr$ defined by 
\[\mu_s(A)=\int_{\gr^0} \#\bigl (s^{-1}(x)\cap A\bigr )d\mu (x)\]
is $\sigma$-finite. 
The analogous statement holds if we replace $s$ by $t$. 
The following conditions on $\mu$ are equivalent. 
\begin{enumerate}[(i)]
\item $\mu_s=\mu_t$, 
\item $i_*\mu_s=\mu_s$, 
\item for every Borel subset $E\subset\gr$ such that 
$s_{\vert E}$ and $t_{\vert E}$ are injective we have 
$\mu(s(E))=\mu(t(E))$. 
\end{enumerate}

\begin{defi}\label{groupoids-definition-measuredgroupoid}
A discrete measurable groupoid $\gr$ together with an invariant
probability measure 
on $\gr^0$, i.e.~satisfying one of (i)-(iii), is 
called a \emph{discrete measured groupoid}. 
\end{defi}

In the sequel $\gr$ will always be a discrete measured groupoid with 
invariant measure $\mu$. The measure on $\gr$ induced by $\mu$ is 
denoted by $\mu_\gr$. 
For a Borel subset $A\subset\gr^0$, 
$\gr_{\vert A}=s^{-1}(A)\cap
t^{-1}(A)$ equipped 
with the normalized measure $\frac{1}{\mu(A)}\mu_{\vert A}$ is 
a discrete measured groupoid, called the \emph{restriction} of 
$\gr$ to $A$. 
The \emph{orbit equivalence relation} on a probability 
$G$-space $(X,\mu)$ defined by  
\[\RC(G\curver X)=\{(x,gx); x\in X, g\in G\}\subset X\times X\] 
is a discrete measured groupoid. The composition is given by 
$(x,y)(y,z)=(x,z)$. 
Another example is given by 
the restriction of the \emph{holonomy groupoid} of a foliation 
to a complete transversal with an invariant measure. 
For a function $\phi:\gr\rightarrow\CC$ and 
$x\in\gr^0$ we put 
\begin{equation*}\begin{aligned}
S(\phi)(x)&=\#\left\{g\in\gr;~\phi(g)\ne
  0,s(g)=x\right\}\in\NN\cup\{\infty\},\\ 
T(\phi)(x)&=\#\left\{g\in\gr;~\phi(g)\ne
  0,t(g)=x\right\}\in\NN\cup\{\infty\}.
\end{aligned}
\end{equation*}
As usual, the set of complex-valued, measurable, essentially bounded 
functions (modulo almost null functions) 
on $\gr$ with respect to $\mu_\gr$ is denoted by 
$L^\infty (\gr,\mu_\gr)$. 
The \emph{groupoid ring} $\CC\gr$ of $\gr$ is defined as 
\[\CC\gr=\left\{\phi\in L^\infty(\gr,\mu_\gr); \text{$S(\phi)$ and $T(\phi)$
    are essentially bounded on $\gr^0$}\right\}.\]
The set $\CC\gr$ is a ring with involution containing 
$\lgrb=L^\infty(\gr^0,\mu)$ 
as a subring. The addition is the pointwise addition 
in $L^\infty (\gr, \mu_\gr)$, the multiplication 
is given by the convolution product 
\[(\phi\eta)(g)=\sum_{\substack{g_1,g_2\in\gr\\g_1g_2=g}}\phi(g_1)\eta(g_2),\quad\phi,\psi\in\CC\gr,~g\in\gr,\]   
and the involution is defined by $(\phi^*)(g)=\overline{\phi(i(g))}$. 
The groupoid ring of the restriction $\CC\gr_{\vert A}$ of 
$\gr$ to $A$, called the \emph{restricted groupoid ring}, 
is canonically isomorphic to $\chi_A\CC\gr\chi_A$. 
Next we explain how $\lgrb$ becomes a left $\CC\gr$-module equipped with 
a $\CC\gr$-epimorphism from $\CC G$ to $\lgrb$. 
The \emph{augmentation homomorphism} is defined by 
$\epsilon:\CC\gr\rightarrow\lgrb$ by 
\[\epsilon:\CC\gr\rightarrow\lgrb,~\epsilon (\phi)(x)=\sum_{g\in
  s^{-1} (x)}\phi (g)\text{ for $x\in\gr^0$}.\]
It becomes a homomorphism of $\CC\gr$-modules when we equip $\lgrb$ with 
the $\CC\gr$-module structure defined below, but it is not a homomorphism 
of rings unless $\gr$ is a group. 
In the language of~\cite{cartan(1999)},  
this means that 
$\CC\gr$ is an augmented ring with the augmentation module $\lgrb$. 
One checks easily that the augmentation homomorphism $\epsilon$ 
induces a $\CC\gr$-module structure 
on $\lgrb$ by  
\[\eta\cdot f=\epsilon (\eta f)\text{ where }\eta\in\CC\gr,~f\in\lgrb.\]
Here $\eta f$ is the product in $\CC\gr$. 

\begin{lem}\label{generators} As a $L^\infty(\gr^0)$-module the
  groupoid ring $\CC\gr$ is generated by the characteristic functions 
$\chi_E$ of Borel subsets $E\subset\gr$ with the property that 
$s_{\vert E}$ and $t_{\vert E}$ are injective. 
\end{lem}

\begin{proof}
Consider $\phi\in\CC\gr$. 
Note that there is some $N>0$ such that 
the preimages of points in $\phi^{-1}(\CC-\{0\})$ under $s$ and $t$ contain 
at most $N$ elements. 
By~\ref{groupoids-definition-kuratowski} 
there is a finite Borel partition $X_i$, $i\in I$,  of 
$\phi^{-1}(\CC-\{0\})$ such that all $s_{\vert X_i}, t_{\vert X_i}$ 
are injective. 
Hence $\phi$ can be written 
as a finite sum $\phi=\sum_{i\in I}\phi_i$, where the support of $\phi_i$ 
lies in $X_i$. Since $s_{\vert X_i}, t_{\vert X_i}$ are 
injective, 
every $\phi_i$ is of the form $f\chi_{X_i}$ (convolution product) 
with $f\in\lgrb$. 
\end{proof}

The groupoid ring $\CC\gr$ of discrete measured groupoid $\gr$ 
lies as a weakly dense involutive
$\CC$-subalgebra in the \emph{von Neumann algebra $\NC(\gr)$ of the
  groupoid $\gr$}. For the construction of $\NC(\gr)$ 
see~\cite{feldman(1977a)},~\cite{feldman(1977b)},
also~\cite{renault(2000)}. 

The von Neumann algebra $\NC(\gr)$ has a finite trace
$\trace_{\NC(\gr)}$ induced by the invariant measure $\mu$. 
For $\phi\in\CC\gr\subset\NC(\gr)$ we have 
\[\trace_{\NC(\gr)}(\phi)=\int_{\gr^0}\phi(g)d\mu(g).\]

Let $R$ be a ring and $G$ be a group. Given a homomorphism 
$c:G\rightarrow\operatorname{Aut}(R)$, $g\mapsto c_g$, we define 
the \emph{crossed product} $R*_cG$ of $R$ with $G$ as the free $R$-module 
with basis $G$. It carries a ring structure that 
is uniquely defined by the rule 
$gr=c_g(r)g$ for $g\in G$, $r\in R$. 
For a probability $G$-space $X$ we denote 
by $\lxbg$ the crossed product $L^\infty (X)*_cG$ obtained by 
the homomorphism $c:G\rightarrow\operatorname{Aut}(L^\infty (X))$ that
maps $g$ to $f\mapsto f\circ l_{g^{-1}}$. 
Here $l_g(x)=gx$ is left translation by $g\in G$. 
The ring homomorphism 
\[\lxbg\rightarrow\CC\RC (G\curver X),~ \sum_{g\in G} f_g\cdot
g\mapsto \bigl ((gx,x)\mapsto f_g(gx)\bigr )\]
is injective and $\phi\in\CC\RC (G\curver X)$ is in the image 
if and only if there is a finite subset $F\subset G$ such that 
$g\not\in F$ implies $\phi (gx, x)=0$ for almost all $x\in X$. 

Note that the map is well defined because the action is 
essentially free. In the sequel we regard $\lxbg$ as a subring of 
$\CC\RC (G\curver X)$. 

\begin{rem} The restriction of the $\CC\RC(G\curver X)$-module 
structure on $L^\infty(X)$ 
to $\lxbg$ is isomorphic to the $\lxbg$-module structure obtained by the 
isomorphism $L^\infty(X)\cong (\lxbg)\otimes_{\CC G}\CC$. 
\end{rem}
\section{Homological Algebra and Dimension
  Theory}\label{sec:some-homol-algebra} 

The objects of study in this section are Tor-groups of the type 
$\tor_n^R(B, M)$ 
where $\BC\subset R\subset\AC$ are ring inclusions, $\AC$, $\BC$
finite von Neumann algebras, $B$ is a $\AC-R$-bimodule and $M$ is 
a $R$-module. Among the questions we deal with are: What happens to 
the $\AC$-dimension of $\tor_n^R(B, M)$ if we replace $M$ by a 
$\dim_\BC$-isomorphic module, $R$ by a $\dim_\BC$-isomorphic subring
$R'\subset R$ or $B$ by a $\dim_\AC$-isomorphic bimodule?

Recall that a ring $R$ is called \emph{semihereditary} if 
every finitely generated submodule of a projective $R$-module is projective. 
A large class of examples for semihereditary rings is given by the 
following theorem. See~\cite[theorem 6.7 on p.~239, p.~288]{lueck(2002)}. 
\begin{thm}\label{groupoids-rings-semihereditary}
Every von Neumann algebra $\NC$ is a semihereditary ring. 
\end{thm}

\begin{thm}[\hbox{\cite[theorem
    1.2]{Lueck(1997a)}}]\label{betti-dim-abeliancategory} 
The category of finitely presented modules 
over a semihereditary ring is abelian. In particular, 
the category of finitely presented modules over a von Neumann algebra 
is abelian.  
\end{thm}

The following theorem is a generalization of 
\cite[theorem 6.29 on p.~253]{lueck(2002)} 
from group von Neumann algebras to arbitrary finite von Neumann algebras. 
A short proof based on the fact that a von Neumann algebra is 
semihereditary can be found in~\cite[theorem
1.48]{sauer(2002)},~\cite[theorem 3.2]{sauer(2003)}.

\begin{thm}\label{groupoids-neumann-flatness}
Any trace-preserving $*$-homomorphism between 
finite von Neumann algebras is a faithfully flat ring extension. 
\end{thm}

\begin{thm}[\hbox{\cite[p.~139-146]{lam(1999)}}]\label{groupoids-rings-flatnesssemihereditary}  
Let $R$ be a semihereditary ring. Then the following holds. 
\begin{enumerate}[(i)]
\item All torsionless $R$-modules are flat.  
\item Any direct product of flat $R$-modules is flat. 
\item Submodules of flat $R$-modules are flat. 
\end{enumerate}
\end{thm}

\begin{lem}\label{groupoids-rings-flatnessgroupoidalgebra}
The groupoid ring $\CC\gr$ of a discrete measured groupoid 
is flat over $\lgrb$. 
\end{lem}

\begin{proof}
There is an inclusion of rings 
$\lgrb\subset\CC\gr\subset\NC(\gr)$ 
where $\NC(\gr)$ is a finite von Neumann algebra whose trace extends 
that of $\lgrb$. By~\ref{groupoids-neumann-flatness}  $\NC(\gr)$ is 
a flat module over $\lgrb$. By the previous theorem $\CC\gr$ is a 
flat $\lgrb$-module. 
\end{proof}

\begin{defi}\label{flipping-bimod}
An $\AC$-$\BC$-bimodule $M$ is called \emph{dimension-compatible} if for every 
$\BC$-module $N$ the following implication holds:
\[\dim_\BC (N)=0\Rightarrow\dim_\AC (M\otimes_\BC N)=0.\]
\end{defi}

We record some easy facts about dimension-compatible bimodules. 

\begin{lem}\label{properties-flippingbimod}\hfill
\begin{enumerate}[(i)]
\item If $M$ is a dimension-compatible $\AC$-$\BC$-bimodule and $N$ is
  a dimension-compatible  
$\BC$-$\mathcal{C}$-bimodule, then $M\otimes_\BC N$ is a dimension-compatible 
$\AC$-$\mathcal{C}$-bimodule. 
\item Quotients and direct summands of dimension-compatible bimodules
  are dimension-compatible. 
\item Let $\BC\subset\AC$ be an inclusion of finite von Neumann algebras. 
Then $\AC$ is a dimension-compatible $\AC$-$\BC$-bimodule. 
\end{enumerate}
\end{lem}

Here the third assertion follows from~\ref{betti-dim-dimpreserving}. 
Next we show that groupoid rings provide examples of dimension-compatible 
bimodules. 

\begin{lem}\label{betti-groupoids-shuffleexample2}
$\CC\gr$ is a dimension-compatible $\lgrb$-$\lgrb$-bimodule. 
\end{lem}

\begin{proof}
Let $M$ be an $\lgrb$-module with $\dim_\lgrb (M)=0$. 
We have to show 
\begin{equation}\label{eq-aim}
\dim_\lgrb\left (\CC\gr\otimes_\lgrb M\right )=0. 
\end{equation}
By the local criterion~\ref{betti-dim-localcriterion} 
equation~\eqref{eq-aim} follows, if  
for a given \mbox{$x\in\CC\gr\otimes_\lgrb M$} a sequence 
of annihilating projections $\chi_{A_i}\in\lgrb$ exists such that  
\begin{equation}\label{eq-generator}
\chi_{A_i}x=0\text{ and }\trace_\lgrb (\chi_{A_i})=\mu(A_i)\rightarrow 1.
\end{equation}
Suppose this is true for a set $S$ of $\lgrb$-generators of 
$\CC\gr\otimes_\lgrb M$. Then \eqref{eq-generator} holds for 
any element in $\CC\gr\otimes_\lgrb M$ by the 
following observation. 
\smallskip\\
If $\chi_{A_i}$ and $\chi_{B_i}$ are annihilating projections 
for the elements $r$ resp.~$s$, whose traces
converge to $1$, then the  
$\chi_{A_i\cap B_i}$ are annihilating projections for $f\cdot r+g\cdot s$, 
$f,g\in\lgrb$, whose traces also converge to $1$. 
\smallskip\\
A set $S$ of $\lgrb$-generators of $\CC\gr\otimes_\lgrb M$ 
is given 
by elements of the form $\chi_E\otimes m$, where $\chi_E$ is 
the characteristic function of a Borel subset 
$E\subset\gr$, such that $s_{\vert E}$ and $t_{\vert E}$ are injective, 
and $m$ is an element in $M$. This is lemma~\ref{generators}. 
Before we prove~\eqref{eq-generator} for the elements in $S$, we show 
that for any Borel subset $A\subset\gr^0$ there is a Borel 
subset $A'\subset\gr^0$ such that 
\begin{equation}\label{eq-prep}
\mu(A')\ge\mu (A)\text{ and }\chi_{A'}\cdot\chi_E=\chi_E\cdot\chi_A.
\end{equation}
We have the identities $\chi_{A'}\chi_E=\chi_{s^{-1}(A')\cap E}$ and 
$\chi_E\chi_A=\chi_{t^{-1}(A)\cap E}$. Put  
\[A'=s\bigl (E\cap t^{-1} (A)\bigr)\cup\bigl (\gr^0-s(E)\bigr).\]
Because $s_{\vert E}$ is injective we get 
\[s^{-1}(A')\cap E=s^{-1}\Bigl (s\bigr (E\cap t^{-1}(A)\bigr )\Bigr)\cap
E=E\cap t^{-1} (A).\]
This yields $\chi_{A'}\cdot\chi_E=\chi_E\cdot\chi_A$. 
The invariance of $\mu$ 
yields 
\begin{equation*}
\begin{split}
\mu (A') &=\mu\Bigl(s\bigl (E\cap t^{-1}(A)\bigr)\Bigr)+\mu
         \bigl(\gr^0-s(E)\bigr )\\
         &=\mu\Bigl (t\bigl (E\cap t^{-1}(A)\bigr)\Bigr)+\mu
         \bigl(\gr^0-t(E)\bigr )\\
         &=\mu\bigl (t(E)\cap A\bigr)+\mu\bigl (\gr^0-t(E)\bigr)\\
         &\ge \mu (A).
\end{split}
\end{equation*}
Now we can prove~\eqref{eq-generator} for 
$x=\chi_E\otimes m\in S$  as follows. 
Because of \mbox{$\dim_\lgrb (M)=0$} 
there are $A_i\subset\gr^0$ with 
$\chi_{A_i}m=0$ and $\mu(A_i)\rightarrow 1$, due to the 
local criterion~\ref{betti-dim-localcriterion}. 
By~\eqref{eq-prep} there are $A_i'\subset\gr^0$ with 
$\chi_{A_i'}\chi_E=\chi_E\chi_{A_i}$ and
$\trace_\lgrb(\chi_{A_i'})=\mu(A_i')\rightarrow 1$. 
Now~\eqref{eq-generator} is obtained from  
$\chi_{A_i'}\left (\chi_E\otimes m\right )=\chi_E\chi_{A_i}\otimes
m=\chi_E\otimes\chi_{A_i}m=0$.
\end{proof}

\begin{lem}\label{betti-groupoids-technical1}\hfill\\
Let $\NC$ be a finite von Neumann algebra, $R$ a ring and $B_1$, $B_2$ 
$\NC$-$R$-bimodules. A bimodule map $B_1\rightarrow B_2$ that is 
a $\dim_\NC$-isomorphism induces $\dim_\NC$-isomorphisms  
\[\tor_\bullet^R (B_1, M)\rightarrow\tor_\bullet^R(B_2, M)\]
for every $R$-module $M$. 
\end{lem}

\begin{proof}
Let $B$ an $\NC$-$R$-bimodule with $\dim_\NC (B)=0$. Let $M$ be an arbitrary 
$R$-module and $P_\bullet$ a free $R$-resolution of $M$. 
Then $\dim_\NC (B\otimes_R P_\bullet)=0$ follows from the 
additivity and cofinality of $\dim_\NC$ (see~\ref{betti-dim-properties}). 
Hence
\[\dim_\NC\left (H_\bullet (B\otimes_R P_\bullet)\right)=\dim_\NC\left
  (\tor_\bullet^R (B, M) \right )=0.\]
In the general case of a $\dim_\NC$-isomorphism 
$\phi:B_1\rightarrow B_2$, we consider the short exact sequences 
\begin{gather*}
0\rightarrow\ker\phi\rightarrow B_1\rightarrow\im\phi\rightarrow 0,\\
0\rightarrow\im\phi\rightarrow B_2\rightarrow\coker\phi\rightarrow 0.
\end{gather*}
$\ker\phi$ and $\coker\phi$ have vanishing dimension. 
We obtain long exact sequences for the Tor-terms: 
\begin{gather*}
\cdots\rightarrow\tor_1^R(B_1,
M)\rightarrow\tor_1^R(\im\phi,
M)\rightarrow\mspace{-15mu}\underbrace{\ker\phi\otimes 
  M}_{=\tor_0^R(\ker\phi,M)}\mspace{-15mu}\rightarrow
B_1\otimes M\rightarrow\im\phi\otimes M\rightarrow 0\\
\cdots\rightarrow\tor_1^R(B_2,
M)\rightarrow\tor_1^R(\coker\phi, M)\rightarrow\im\phi\otimes M\rightarrow
B_2\otimes M\rightarrow\mspace{-15mu}\underbrace{\coker\phi\otimes
  M}_{=\tor_0^R(\coker\phi,M)}\mspace{-15mu}\rightarrow 0.
\end{gather*}
We already know 
$\dim_\NC\left (\tor_\bullet^R(\ker\phi,M)\right )=0$ and 
\mbox{$\dim_\NC\left (\tor_\bullet^R(\coker\phi,M)\right )=0$}, hence 
\begin{gather*}
\tor_\bullet^R(B_1, N)\rightarrow\tor_\bullet^R (\im\phi, N),\\
\tor_\bullet^R(\im\phi, N)\rightarrow\tor_\bullet^R (B_2, N)
\end{gather*}
are $\dim_\NC$-isomorphisms, and so is their composition. 
\end{proof}

\begin{lem}\label{betti-groupoids-technical2}\hfill\\
Let $\BC\subset R\subset\AC$ be an inclusion of rings where 
$\AC$, $\BC$ are finite von Neumann algebras. Let $B$ be an
$\AC$-$R$-bimodule. We assume the following.  
\begin{enumerate}[(i)]
\item $R$ is dimension-compatible as a $\BC$-$\BC$-bimodule. 
\item $B$ is dimension-compatible  as an $\AC$-$\BC$-bimodule. 
\item $B$ is flat as a right $\BC$-module. 
\end{enumerate}
Then every $R$-homomorphism $M\rightarrow N$, which 
is a $\dim_\BC$-isomorphism, induces 
a $\dim_\AC$-iso\-morphism
\[\tor_\bullet^R (B, M)\rightarrow\tor_\bullet^R (B, N).\]
\end{lem}

\begin{proof}
First we get the following flatness properties. 
\begin{enumerate}[$\bullet$]
\item The $\AC$-$R$-bimodule $B\otimes_\BC R$ is flat as a right $R$-module 
because $B$ is flat as a right $\BC$-module. 
\item $R$ is a $\BC$-submodule of the flat $\BC$-module $\AC$ 
(\ref{groupoids-neumann-flatness}). Hence $R$ is flat as a right 
$\BC$-module 
by~\ref{groupoids-rings-flatnesssemihereditary} and 
\ref{groupoids-rings-semihereditary}. 
\item Therefore $B\otimes_\BC R$ is also flat as a right $\BC$-module. 
\end{enumerate}
Multiplication yields a 
surjective $\AC$-$R$-bimodule homomorphism 
\[m: B\otimes_\BC R\rightarrow B.\]
$B\otimes_\BC R$ is dimension-compatible (as an $\AC$-$\BC$-bimodule) because 
$B$ and $R$ are dimension-compatible. 
The map $m$ splits as an $\AC$-$\BC$-homomorphism by the map 
$B\rightarrow B\otimes_\BC R$, $b\mapsto b\otimes 1$. Hence, as
an $\AC$-$\BC$-bimodule, $\ker m$ is a direct summand in 
$B\otimes_\BC R$ and therefore also dimension-compatible. 
We record the properties of $\ker m$: 
\begin{enumerate}[$\bullet$]
\item $\ker m$ is an $\AC$-$R$-bimodule. 
\item $\ker m$ is dimension-compatible as an $\AC$-$\BC$-bimodule.  
\item $\ker m$ is flat as a right $\BC$-module because it is the submodule 
of a flat $\BC$-module 
(\ref{groupoids-rings-flatnesssemihereditary},
\ref{groupoids-rings-semihereditary}). 
\end{enumerate}
Notice that $\ker m$ satisfies all properties imposed 
on $B$. Let $M$ be an $R$-module with $\dim_\BC (M)=0$. 
So we have $\dim_\AC\left (\ker m\otimes_\BC M\right )=0$ and 
hence for its quotient
\[\dim_\AC\left (\ker m\otimes_R M\right )=0.\]
The short exact sequence 
$0\rightarrow\ker m\rightarrow B\otimes_\BC R\rightarrow B\rightarrow 0$ 
induces a long exact sequence of Tor-terms 
\begin{multline*}
\ldots\rightarrow 0\rightarrow\tor_2^R (B, M)\rightarrow
\tor_1^R (\ker m, M)\rightarrow\underbrace{\tor_1^R(B\otimes_\BC R,M)}_{=0}
\rightarrow\\\rightarrow\tor_1^R (B, M)\rightarrow\ker m\otimes_R M\rightarrow
\left (B\otimes_\BC R\otimes_R M\right )\rightarrow B\otimes_R M\rightarrow 0,
\end{multline*}
where the zero terms are due to $B\otimes_\BC R$ being $R$-flat. 
We obtain  
\begin{equation}\label{eq-repeat}
\begin{gathered}
\dim_\AC\left (B\otimes_R M\right )=\dim_\AC\left (\tor_1^R (B, M)\right
)=0,\\
\tor_{i+1}^R (B, M)\cong\tor_i^R(\ker m, M)\quad i\ge 1.
\end{gathered}
\end{equation}
Now we apply~\eqref{eq-repeat} to 
$\ker m$ instead of $B$ and get 
\[\dim_\AC\left (\tor_1^R (\ker m, M)\right )=0,~\dim_\AC\left (\tor_2^R (B,
  M)\right )=0.\]  
Repeating this ("dimension shifting") yields 
\[\dim_\AC\left (\tor_i^R (B, M)\right )=0\text{ for $i\ge 0$.}\] 
Deducing the general case of a 
$\dim_\BC$-isomorphism $\phi:M\rightarrow N$ from the case $\dim_\BC (M)=0$ 
uses exactly the same method as in the proof of 
\ref{betti-groupoids-technical1}. 
\end{proof}

\begin{thm}\label{betti-groupoids-smallerring}\hfill\\
Let $\BC\subset R_1\subset R_2\subset\AC$ be an 
inclusion of rings where $\AC$, $\BC$ are finite von Neumann algebras.  
We assume the following. 
\begin{enumerate}[(i)]
\item $R_2$ is dimension-compatible as a $\BC$-$\BC$-bimodule. 
\item The inclusion $R_1\subset R_2$ is a $\dim_\BC$-isomorphism. 
\end{enumerate} 
Then 
\[\dim_\AC\left (\tor_\bullet^{R_1} (\AC, M)\right )=\dim_\AC\left
  (\tor_\bullet^{R_2} (\AC, M)\right )\]
holds for every $R_2$-module $M$. 
\end{thm}

\begin{proof}
By~\ref{betti-groupoids-technical1} the induced map 
\begin{equation}\label{eq-dimiso}
\tor_i^{R_1} (R_2, M)\longleftarrow\tor_i^{R_1} (R_1, M)=
\begin{cases}
           M&\text{ if $i=0$}\\
           0&\text{ if $i>0$}
         \end{cases}
\end{equation}
is a $\dim_\BC$-isomorphism. 
Let $P_\bullet\rightarrow M$ be a projective $R_1$-resolution of $M$. 
The K\"unneth spectral sequence, applied to $\AC$ and the complex
$R_2\otimes_{R_1} P_\bullet$ (see~\cite[theorem 5.6.4 on
p.~143]{Weibel(1994)}), has the $E^2$-term 
\[E_{pq}^2=\tor_p^{R_2} \Bigl(\AC, H_q (R_2\otimes_{R_1} P_\bullet)\Bigr)=
 \tor_p^{R_2} \Bigl (\AC,\tor_q^{R_1}(R_2,M)\Bigr)\]
and converges to 
\[H_{p+q} \Bigl (\AC\otimes_{R_2}\bigl (R_2\otimes_{R_1}
P_\bullet\bigr )\Bigr )=\tor_{p+q}^{R_1} (\AC,M).\]
Now we can apply 
the preceding lemma~\ref{betti-groupoids-technical2} 
to the inclusion $\BC\subset R_2\subset\AC$. Here recall that $\AC$ is flat as 
a right $\BC$-module (\ref{groupoids-neumann-flatness}) and
dimension-compatible as an $\AC$-$\BC$-bimodule. So we obtain 
\[\dim_\AC \left (E_{pq}^2\right )=
\begin{cases} 
\dim_\AC\Bigl (\tor_p^{R_2} (\AC, M)\Bigr )  & \text{ if $q=0$}\\
0                              & \text{ if $q>0$.}
\end{cases}\]
So the spectral sequence collapses up to dimension, and  
additivity (\ref{betti-dim-properties}) yields
\[\dim_\AC\left (\tor_p^{R_1}(\AC, M)\right )=
\dim_\AC\left ( E^\infty_{p0}\right )=\dim_\AC
\left (E^2_{p0}\right )=\dim_\AC\left (\tor_p^{R_2} (\AC, M)\right ).\]
\end{proof}

Actually, the theorem above provides more than just an equality of the
dimensions. There is a natural zig-zag of $\dim_\AC$-isomorphisms 
between $\tor_\bullet^{R_1} (\AC, M)$ and 
$\tor_\bullet^{R_2} (\AC, M)$. 

\begin{thm}\label{betti-groupoids-reducedtor}\hfill\\
Let $\BC\subset R\subset\AC$ be an inclusion of rings, where 
$\AC$, $\BC$ are finite von Neumann algebras. Let $p\in\BC$ be a projection. 
We assume the following. 
\begin{enumerate}[(i)]
\item $R$ is dimension-compatible as a $\BC$-$\BC$-bimodule. 
\item $p$ is $\dim_\BC$-full in $R$. 
\end{enumerate}
Then the equality 
\begin{gather*}
\dim_{p\AC p}\Bigl (p\tor_\bullet^R (\AC, M)\Bigr )=\dim_{p\AC p}\Bigl
(\tor_\bullet^{pRp}(p\AC p, pM)\Bigr )
\end{gather*}
holds for every $R$-module $M$. 
\end{thm}

\begin{proof}
We have $p\tor_i^{pRp} (Rp, pM)\cong\tor_i^{pRp} (pRp,pM)=0$ for
$i>0$. Theorem~\ref{almost-morita} implies 
$\dim_\BC\bigl (\tor_i^{pRp}(Rp, pM)\bigr )=0$, $i>0$. 
By the same theorem the multiplication map $m$ 
\[\tor_0^{pRp}(Rp,pM)=Rp\otimes_{pRp} pM\longrightarrow M\] 
is a $\dim_\BC$-isomorphism. 
Now let $P_\bullet\rightarrow pM$ be a $pRp$-projective resolution of 
$pM$. 
The K\"unneth spectral sequence applied to $\AC$ and the complex
$Rp\otimes_{pRp} P_\bullet$ (\cite[theorem 5.6.4 on
p.~143]{Weibel(1994)}) 
has the  $E^2$-term
\[E_{ij}^2=\tor_i^R\Bigl (\AC, H_j(Rp\otimes_{pRp}
P_\bullet)\Bigr )=\tor_i^R\Bigl (\AC,\tor_j^{pRp}(Rp, pM)\Bigr)\]
and converges to 
\[H_{i+j} \bigl(\underbrace{\AC\otimes_R (Rp\otimes_{pRp}
P_\bullet)}_{=\AC p\otimes_{pRp}P_\bullet}\bigr)=\tor_{i+j}^{pRp}\bigl (\AC p,
pM\bigr).\]
We know $\tor_\bullet^{pRp} (Rp, pM)$ up to 
$\dim_\BC$-isomorphism. An application of 
Lemma~\ref{betti-groupoids-technical2} to that module 
yields  
\[\dim_\AC \left (E^2_{ij}\right )=\begin{cases} 
   \dim_\AC\Bigl (\tor_i^R (\AC, M)\Bigr ) &\text{ if $j= 0$}\\
   0                                     &\text{ if $j>0$.}
 \end{cases}\]
The spectral sequence collapses up to dimension. This implies  
\[\dim_\AC\left (\tor_i^{pRp} (\AC p, pM)\right )=\dim_\AC \left
  (E^\infty_{i0}\right )=\dim_\AC\left ( E^2_{i0}\right )=\dim_\AC\left
  (\tor_i^{R} (\AC, M)\right ).\]
It follows from~\ref{betti-dim-localcriterion} that $\AC
  p\AC\subset\AC$ is $\dim_\AC$-surjective since $RpR\subset R$ is
  $\dim_\BC$-surjective. Hence $p$ is $\dim_\AC$-full in $\AC$ and the
  claim is obtained from~\ref{betti-dim-projectiondim}. 
\end{proof}

\section{L$^2$-Betti Numbers of Discrete Measured
  Groupoids}\label{sec:l2-betti-numbers}

\begin{defi}\label{betti-definition}
Let $\gr$ be a discrete measured groupoid. Its \emph{n-th $L^2$-Betti 
number} $\betti_n(\gr)$ is defined as 
\[\betti_n(\gr)=\dim_{\NC(\gr)}\left (\tor_n^{\CC\gr}\left (\NC (\gr),L^\infty
(\gr^0)\right )\right ).\]
\end{defi}

\begin{lem}\label{betti-restriction0}
Let $\gr$ be a discrete measured groupoid, and let $A\subset\gr^0$ 
be a Borel subset such that $t(s^{-1}(A))$ has full measure in 
$\gr^0$. Then the characteristic function $\chi_A\in\lgrb\subset\CC\gr$ of $A$ 
is a $\dim_\lgrb$-full idempotent in $\CC\gr$. 
\end{lem}

\begin{proof}
We have to show that the inclusion 
$\CC\gr\chi_A\CC\gr\subset\CC\gr$ is a
$\dim_{L^\infty(\gr^0)}$-isomorphism.  
Let $s^{-1}(A)=\bigcup_{n\in\NN} E_n$ be a partition into Borel 
sets $E_n$ with the property that $s_{\vert E_n}$ is injective. 
The partition exists by 
theorem~\ref{groupoids-definition-kuratowski}. 
We have 
\[\chi_{i(E_n)}\cdot\chi_A\cdot\chi_{E_n}=\chi_{i(E_n)}\chi_{s^{-1}(A)\cap
  E_n}=\chi_{t(s^{-1}(A)\cap E_n)}.\]
The right equal sign is due to the injectivity of $s_{\vert E_n}$. 
Hence 
$\sum_{n=1}^N\chi_{t(s^{-1}(A)\cap E_n)}\in\CC\gr\chi_A\CC\gr$. 
So we get  
\[\chi_{t(s^{-1}(A)\cap\bigcup_{n=1}^N E_n)}=f\cdot\left
  (\sum_{n=1}^N\chi_{t(s^{-1}(A)\cap E_n)}\right
  )\in\CC\gr\cdot\chi_A\cdot\CC\gr\]
for a suitable $f\in L^\infty(\gr^0)$. 
This implies $\chi_{t(s^{-1}(A)\cap\bigcup_{n=1}^N E_n)}\cdot\phi=0$ 
for 
every element $\phi$ in the quotient $\CC\gr/\CC\gr\chi_A\CC\gr$. 
Because $t(s^{-1}(A))$ has full measure, we get 
\[\lim_{N\rightarrow\infty}\mu\Bigl (t\bigl (s^{-1}(A)\cap\bigcup_{n=1}^N
  E_n\bigr)\Bigr)= \mu\bigl (\gr^0\bigr)=1.\]
So $\CC\gr\chi_A\CC\gr\subset\CC\gr$ is a
$\dim_{L^\infty(\gr^0)}$-isomorphism
by~\ref{betti-dim-localcriterion}. 
\end{proof}

\begin{thm}\label{betti-restriction}
Let $\gr$ be a discrete measured groupoid, and let $A\subset\gr^0$ 
be a Borel subset such that $t(s^{-1}(A))$ has full measure in 
$\gr^0$. Then 
\[\betti_n(\gr)=\mu (A)\cdot\betti_n (\gr_{\vert A}).\]
\end{thm}

\begin{proof} By~\ref{betti-restriction0} 
and~\ref{betti-dim-projectiondim} we get 
\[\betti_n (\gr)=\mu(A)\cdot\dim_{\chi_A\NC(\gr)\chi_A}\Bigl
(\chi_A\tor_n^{\CC\gr}\bigl (\NC (\gr),L^\infty (\gr^0)\bigr )\Bigr).\]
By~\ref{betti-groupoids-shuffleexample2} the groupoid ring $\CC\gr$ is
  dimension-compatible as a $\lgrb$-$\lgrb$-bimodule whence by 
\ref{betti-groupoids-reducedtor} we obtain that $\betti_n(\gr)$ equals 
\[\mu(A)\cdot\dim_{\chi_A\NC(\gr)\chi_A}\Bigl
  (\tor_n^{\chi_A\CC\gr\chi_A} \bigl 
  (\chi_A\NC(\gr)\chi_A,\chi_A\lgrb\bigr )\Bigr )=
\mu(A)\cdot\betti_n \bigl (\gr_{\vert A}\bigr ).\]
For the right equal sign note that 
$\NC (\gr_{\vert A})=\chi_A\NC(\gr)\chi_A,~\CC\gr_{\vert
  A}=\chi_A\CC\gr\chi_A$ and $\chi_AL^\infty (\gr^0)=L^\infty(A)$. 
\end{proof}

\begin{lem}\label{crossedintogroupoid}
Let $X$ be a probability $G$-space and be $\RC$ the orbit equivalence
relation on $X$. Then the 
inclusion $\lxbg\subset\CC\RC$ is a $\dim_\lgrb$-isomorphism. 
\end{lem}

\begin{proof}
We apply the local criterion
(\ref{betti-dim-localcriterion}) to show that 
the quotient $\CC\RC/\lxbg$ has dimension zero. 
Let $\phi\in\CC\RC$ and $[\phi]$ be its image in the quotient. 
Choose an enumeration 
$G=\{g_1,g_2,\ldots\}$. Define Borel 
subsets $X_n\subset X$ by 
$X_n=\{x\in X;~\phi (g_ix,x)=0\text{ for all $i>n$}\}$, and let 
$\chi_n$ be the characteristic function of $X_n$. 
Then $\chi_n\phi\in\lxbg$ holds, hence 
$\chi_n[\phi]=0\in\CC\RC/\lxbg$. Because of $\mu(X_n)\rightarrow\mu
(X)=1$ the local criterion yields $\dim_\lxb\left (\CC\RC/\lxbg\right )=0$. 
\end{proof}

\begin{thm}\label{betti-mainthm}
Let $X$ be a probability $G$-space. 
Then the $L^2$-Betti numbers of $G$ and the orbit equivalence relation 
$\RC$ on $X$ coincide. 
\[\betti_n (G)=\betti_n\left(\RC\right )\]
\end{thm}

\begin{proof} 
The crossed product ring $\lxbg$ 
is flat as a right $\CC G$-module  
because of the equality 
$\left (\lxbg\right )\otimes_{\CC G} M=L^\infty (X)\otimes_\CC M$. 
Hence we obtain 
\begin{align*} 
\betti_n (G) &=\dim_{\NC(\RC)}\left (\NC
 (\RC)\otimes_{\NC (G)}\tor_n^{\CC G}\left (\NC (G),\CC\right )\right
  )\qquad\text{by (\ref{betti-dim-dimpreserving})}\\ 
&=\dim_{\NC(\RC)}\left (\tor_n^{\CC G}\left (\NC (\RC),\CC\right )\right
  )\qquad\text{by (\ref{groupoids-neumann-flatness})}\\ 
 &=\dim_{\NC(\RC)}\left (\tor_n^\lxbg\left (\NC (\RC),\left (\lxbg\right
  )\otimes_{\CC G}\CC\right )\right ) \\
 &=\dim_{\NC(\RC)}\left (\tor_n^\lxbg\left (\NC (\RC),L^\infty
  (X)\right)\right ).
\end{align*}
Now we will apply theorem~\ref{betti-groupoids-smallerring} 
to the ring inclusions 
\[L^\infty (X)\subset\lxbg\subset\CC\RC\subset\NC (\RC).\]
As an $\lxb$-$\lxb$-bimodule, $\CC\RC$ is dimension-compatible 
by~\ref{betti-groupoids-shuffleexample2}. The inclusion 
$\lxbg\subset\CC\RC$ is a $\dim_\lgrb$-isomorphism by
\ref{crossedintogroupoid}. Hence by~\ref{betti-groupoids-smallerring} 
we obtain 
\begin{equation*}
\begin{split}
\dim_{\NC(\RC)}\left (\tor_n^\lxbg\left (\NC (\RC),L^\infty (X)\right
  )\right) &= 
\dim_{\NC(\RC)}\left (\tor_n^{\CC\RC}\left (\NC (\RC),L^\infty
  (X)\right )\right )\\
&=\betti_n(\RC).
\end{split}
\end{equation*}
\end{proof}

\bibliographystyle{alpha}
\bibliography{promotion}
\end{document}